\documentclass[12pt]{amsart}
\usepackage{epsf}

\begin{document}

\title{Manifolds not containing Gompf nuclei}
\author{Andr\'{a}s I. Stipsicz}
\thanks{The author was supported by the Magyary Zolt\'an
Foundation and by OTKA.
Research at MSRI is supported in part by NSF grant DMS-9022140.}
\address{Department of Analysis,\\
ELTE TTK, Budapest, Hungary}
\email{stipsicz@cs.elte.hu}

\newcommand{\cp}{{\Bbb {CP}}}
\newcommand{\bfz}{{\Bbb {Z}}}
\newtheorem{fact}{Fact}[section]
\newtheorem{lemma}[fact]{Lemma}
\newtheorem{theorem}[fact]{Theorem}
\newtheorem{definition}[fact]{Definition}
\newtheorem{remark}[fact]{Remark}
\newtheorem{proposition}[fact]{Proposition}

\newenvironment{prooff}{\medskip \par \noindent {\it Proof}\ }{\hfill
$\square$ \medskip \par}
	\def\sqr#1#2{{\vcenter{\hrule height.#2pt
   		\hbox{\vrule width.#2pt height#1pt \kern#1pt
      		\vrule width.#2pt}\hrule height.#2pt}}}
	\def\square{\mathchoice\sqr67\sqr67\sqr{2.1}6\sqr{1.5}6}
\def\pf#1{\medskip \par \noindent {\it #1.}\ }
\def\endpf{\hfill $\square$ \medskip \par}
\def\demo#1{\medskip \par \noindent {\it #1.}\ }
\def\enddemo{\medskip \par}
\def\qed{~\hfill$\square$}

\font\sevenrm=cmr7
\font\tenmsb=msbm10  at 11pt
\font\sevenmsb=msbm7 at 9pt
\font\fivemsb=msbm5 at 7pt
\newfam\msbfam
\textfont\msbfam=\tenmsb
\scriptfont\msbfam=\sevenmsb
\scriptscriptfont\msbfam=\fivemsb
\def\Bbb#1{\fam\msbfam\relax#1}
	\def\complex{{\Bbb C}}
	\def\IF{{\Bbb F}}
	\def\que{{\Bbb Q}}
	\def\pee{{\Bbb P}}
	\def\real{{\Bbb R}}
	\def\nat{{\Bbb N}}
	\def\zed{{\Bbb Z}}
\def\C{{\cal C}}
\def\O{{\cal O}}
\def\cS{{\cal S}}
\def\simtimes{\,\tilde\times\,}	
\def\chix{{\raise.5ex\hbox{$\chi$}}}
\def\frac#1#2{{\textstyle{#1\over#2}}}
\def\cpn{{\complex {\Bbb {P}}^n}}
\def\cpone{{\complex {\Bbb {P}}^1}}
\def\cptwo{{\complex {\Bbb {P}}^2}}
\def\ocptwo{\overline{\complex {\Bbb {P}}}^2} 
\def\blk{$\underline{\qquad}$} 
\def\ep{\varepsilon}
\def\ov{\overline}

\begin{abstract}
In this note we show that there are 4-manifolds not containing
Gompf nucleus $N_2$; in this way we answer 
Problem 4.98 of Kirby's problem list (see \cite{K}) in the 
negative.
\end{abstract}
\maketitle

\section{Introduction}
\label{elso}

This note is devoted to answer a question in Kirby's problem
list \cite{K} asking whether every simply connected 
smooth 4-manifold with $b^+\geq 3$ contains a Gompf nucleus
$N_2$ (Problem 4.98 in \cite{K}). We prove that by doing 
logarithmic transformations on three linearly independent
 tori in the {\em K3\/}-surface we get a 4-manifold not contaning
$N_2$. To make our statements more precise we need a few definitions.

\begin{figure}
\centerline{		
\epsfbox{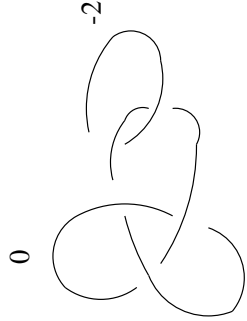}}
\caption{Kirby picture for $N_2$}		
\label{egyes}
\end{figure}

The hypersurface $X=\lbrace [z_0:z_1:z_2:z_3]\in \cp ^3\ \vert \ 
\sum _{i=0} ^3 z_i ^4=0\rbrace \subset \cp ^3$ 
is a simply connected, smooth 4-manifold with $c_1(X)=0$, hence it is
a {\em K3-surface\/}. It is known that all simply connected complex surfaces
with vanishing first Chern class are diffeomorphic, consequently
from the differential topological point of view $X$ is {\em the K3\/}-surface.
The complex surface $X$ admits a holomorphic fibration $\pi \colon
X\to \cp ^1$ such that the generic fiber is a smooth elliptic curve
--- a 2-dimensional torus. Such fibrations are called {\em elliptic 
fibrations\/}. It can be assumed that $\pi $ has a singular fiber homeomorphic
to the 2-dimensional sphere $S^2$ --- such fibers are called 
{\em cusp fibers\/}.
The fibration also has a section $\sigma \colon \cp ^1\to X$, the image
of which is a sphere $S\subset X$ with square $[S]^2=-2$. The 
{\em Gompf nucleus\/} $N_2$ is by definition the tubular neighborhood
of the union of a cusp fiber and a section in $X$. The manifold $N_2$ 
admits a handle decomposition with one 0-handle and two 2-handles, where
these two 2-handles are attached to $S^3=\partial (0-handle)$ according to
the Kirby diagram shown in Figure~\ref{egyes}.

It can be shown that the {\em K3\/}-surface $X$ contains three disjoint
copies of $N_2$. If a 4-manifold $M^4$ contains a 2-dimensional
torus $T$ with square 0, then we can perform a {\em logarithmic 
transformation\/} on $T$: deleting the tubular neighborhood of $T$
(which is diffeomorphic to $D^2\times T^2$) and regluing it via a 
diffeomorphism $\varphi \colon \partial (M\setminus D^2\times T^2 )\to 
\partial (D^2\times T^2)$ we get a new manifold $M_{\varphi}$.
It turnes out that if $T$ is the fiber 
in a Gompf nucleus $N_2\subset M$, then the
diffeomorphism type of $M_{\varphi}$ will depend only on one nonnegative 
number $p$ associated to $\varphi$. This number is called the
{\em multiplicity\/} of the logarithmic transformation.
For more about elliptic surfaces, nuclei and logarithmic transformations
see \cite{G}, \cite{FS1} and \cite{GS}.

Now perform logarithmic transformations of multiplicity 2 on the three
tori contained by the three disjoint nuclei in the {\em K3\/}-surface
$X$. The resulting manifold will be denoted by $X_{2,2,2}$.

\begin{theorem}
The 4-manifold $X_{2,2,2}$ does not contain a Gompf nucleus $N_2$.
\label{elsotetel}
\end{theorem}

\begin{remark}
{\rm One of the most interesting questions in 4-manifold theory is whether
all 4-manifolds are of simple type or not. (For the definition
of simple type see Section~\ref{masodik}.) It is known that if $M^4$
contains a homologically essential torus with square 0, then
$M$ is of simple type. There is no example of a 4-manifold with $b^+\geq 3$
and not containing a torus with square 0. One can ask whether every 
4-manifold contains a cusp neeighborhood (a tubular neighborhood
of a cusp fiber) --- or even a Gompf nucleus $N_2$. The above theorem
shows an example of a manifold which contains no $N_2$; $X_{2,2,2}$ is still
of simple type, however.}
\end{remark}
 
\begin{figure}[b]
\centerline{		
\epsfbox{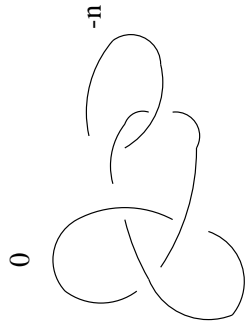}}
\caption{Kirby picture for $N_n$}		
\label{kettes}
\end{figure}

One can modify the question by trying to find more general nuclei
$N_n$ in 4-manifolds. The manifold $N_n$ is described by the 
Kirby diagram shown in Figure~\ref{kettes} --- it can be defined alternatively
as the tubular neighborhhod of the union of a cups fiber and a section
in a simply connected elliptic surface (admitting a 
section) with Euler characteristic 
$12n$ (see \cite{G}).

\begin{theorem}
For every $n$ there is a manifold $Y_n$ such that $Y_n$ does not contain 
$N_n$.
\label{masodiktetel}
\end{theorem}

A related question would be to find a 4-manifold $M$ with the property that
$M$ does not contain $N_2(p,q)$; here $N_2(p,q)$ denotes the 4-manifold
with boundary we get by performing two logarithmic transformations 
(of multiplicity $p$ and $q$) along the fiber of the Gompf nucleus
$N_2$. Using a recent construction of Fintushel and Stern \cite{FS3}
such $M$ can be found, see \cite{SSz}.

{\em Acknowledgement}: The author would like to thank MSRI
for their hospitality and support, and Zolt\'an Szab\'o
for many helpful conversations.

\section{Basic classes}
\label{masodik}

In studying differential topological properties of smooth 
4-manifolds, {\em Seiberg-Witten basic classes\/} 
turn out to be very powerful tools. For the definition
of these objects see \cite{A}, \cite{M} or \cite{GS},
here we restrict ourselves to a very short outline.

Assume that $M^4$ is a smooth, oriented, simply connected, closed
4-manifold with $b^+_M \geq 3$. A cohomology class $K\in H^2 (M;\bfz )$
with $K\equiv w_2(M)\ (mod\ 2)$ uniquely determines a spin$^c$ structure
on $M$, and for such a structure a certain pair of 
partial differential equations
(the so-called Seiberg-Witten equations)
can be described --- involving a choice of metric on $M$, a coupled 
Dirac operator and a perturbation 2-form. By a delicate "counting argument"
of the solutions of these equations a number $SW_M(K)$
can be associated to the cohomology class $K$. It turnes out that this number
(up to sign) is a smooth invariant of the manifold $M$, more precisely

\begin{theorem}
If $f\colon M'\to M$ is an orientation preserving diffeomorphism
then $SW_M(K)=\pm SW _{M'}(f^{\ast}K)$. Moreover, for a fixed 4-manifold
$M$ there are only finitely many classes $K$ with $SW_M(K)\neq 0$,
and $SW_M(-K)=\pm SW_M(K)$.\qed
\end{theorem}

\begin{definition}
{\rm 
The cohomology class $K\in H^2(M;\bfz )$ is called a 
{\em Seiberg-Witten basic class\/} if $SW_M(K)\neq 0$. The 4-manifold
$M$ is of simple type if $SW_M(K)\neq 0$ implies that 
$K^2=3\sigma (M)+\chi (M)$; here $\sigma (M)$ and $\chi (M)$ stand for the
signature and Euler characteristic of $M$.}
\end{definition}

The most important relation between the smooth topology of a 4-manifold $M$
and its basic classes $\lbrace K_i \ \vert \ i=1,...,n\rbrace $ is
shown by the {\em generalized adjunction formula\/}:

\begin{theorem} {\rm {(Kronheimer-Mrowka)}}
If $\Sigma ^2\subset M^4$ is a smooth, connected 2-dimen\-sional surface of genus
$g(\Sigma )$, $[\Sigma ]\neq 0$  
and $[\Sigma ]^2\geq 0$, then for every basic class
$K$ we have 
$$ 2g(\Sigma )-2 \geq [\Sigma ]^2+\vert K([\Sigma ])\vert .$$ \qed
\end{theorem}

\section{Proofs of Theorems~\ref{elsotetel} and \ref{masodiktetel}}
\label{harmadik}

Assume that $T_1,T_2$ and $T_3$ are three tori in the 
{\em K3\/}-surface lying in three disjoint Gompf 
nuclei. Perform logarithmic transformations of multiplicity 2
on each $T_i$. The basic classes of the resulting manifold
$X_{2,2,2}$ are determined by Fintushel and Stern \cite{FS2}.

\begin{proposition}
The basic classes of $X_{2,2,2}$ are the Poincar\'e duals
of the homology classes $\pm [{{T_1}\over 2}] \pm  [{{T_2}\over 2}]
\pm  [{{T_3}\over 2}]$.\qed
\end{proposition}

\begin{prooff} {\em of Theorem~\ref{elsotetel}:\/}
Assume that $N_2\subset X_{2,2,2}$. The homology class of the 
fiber and the section in $N_2$ will be denoted by $f$ and
$s$ respectively; note that $f$ and $g=f+s$ have square 0 and can be represented by tori.
Consequently (by the generalized adjunction formula)
a basic class $K$ of $X_{2,2,2}$ evaluates trivially on $f$ and $g$.
Now $f\cdot  ([{{T_i}\over 2}]+ [{{T_j}\over 2}]+ [{{T_k}\over 2}])=0$
and $f\cdot (-[{{T_i}\over 2}]+ [{{T_j}\over 2}]+ [{{T_k}\over 2}])0$
($\lbrace i,j,k\rbrace =\lbrace 1,2,3\rbrace$)
implies that $f\cdot [{{T_i}\over 2}]=0$, similarly $g\cdot [{{T_i}\over 2}]
=0$ for $i=1,2,3$. Since the complement of the three disjoint nuclei
in $X$ have negative definite intersection form, the above equalities
imply that $f=\sum _{i=1} ^3 \alpha _i [T_i]$ and similarly
$g=\sum _{i=1} ^3 \beta _j [T_j]$. 
These two latter equations, however,
give a contradiction, since $f\cdot g=1$ but
$(\sum _{i=1} ^3 \alpha _i [T_i])\cdot (\sum _{i=1} ^3 \beta _j [T_j])=0$. 
Consequently $N_2$ does not embed in $X_{2,2,2}$.
\end{prooff}

\begin{prooff} {\em of Theorem~\ref{masodiktetel}:\/}
Perform logarithmic transsformations of multiplicity $2n$ on $T_1,T_2$ and 
$T_3$ as above; the resulting 4-manifold $X_{2n,2n,2n}$ will be denoted by
$Y_n$. The basic classes of $Y_n$ are determined in \cite{FS2}:
these are the Poincar\'e duals 
of the homology classes
$\alpha [{{T_1}\over {2n}}]+\beta [{{T_2}\over {2n}}]+
\gamma [{{T_3}\over {2n}}]$ where $\alpha ,\beta ,\gamma \in 
\lbrace \pm (2j-1) \ \vert \ j=1,...,n\rbrace $.
Assume now that $N_n\subset Y_n$. As before, $f$ and $s$ will 
denote the homology classes of the fiber and the section in $N_n$
respectively. Note that $f^2=0$ and $s^2=-n$. The class $f$ can be 
represented by a torus, while $g=s+nf$ can be represented
by a surface of genus $n$; moreover $g^2=n$. The same argument 
as in the proof of Theorem~\ref{elsotetel}
gives that $f=\sum _{i=1} ^3 \alpha _i [T_i]$. To show that every 
basic class of $Y_n$ evaluates trivially on $g$ involves a little
trick: take the basic
classes $\pm [{{T_1}\over 2}]\pm [{{T_2}\over 2}]\pm [{{T_3}\over 2}]$
and evaluate them on $g$. If one of them, say $K$, evaluates nontrivially
on $g$, then for $L=(2n-1)K$ (which is also a basic class)
we have $\vert L(g)\vert \geq 2n-1$, but then $L$ and $g$ would violate
the generalized adjunction formula. Consequently
$(\pm [{{T_1}\over 2}]\pm [{{T_2}\over 2}]\pm [{{T_3}\over 2}])\cdot g=0$
implying that $g\cdot [T_i]=0$ $(i=1,2,3)$.
These latter equations show that $g=\sum _{j=1} ^3 \beta _j [T_j]$, and we 
have the same contradiction as before.
\end{prooff} 

\begin{remark}
{\rm 
Note that the proof given above also shows that no 
$N_k$ with $k\leq n$ embeds in $Y_n=X_{2n,2n,2n}$.}
\end{remark}

\end{document}